\documentclass[final,leqno]{siamltex704}
\usepackage{amsmath}
\usepackage{amsfonts}
\usepackage{latexsym}
\usepackage{amssymb}

\renewcommand{\ldots}{\dotsc}
\newcommand{\R}{\mathbb{R}}

\def \d {\displaystyle}

\def \sysd {\left\{\begin{array}{l}}
\def \sysf {\end{array}\right.}
\def \eqd {\begin{equation}}
\def \eqf {\end{equation}}
\def \f {\end{document}}

\def\S{\scriptstyle}

\newtheorem{remark}{\textit{Remark}}

\title{Global Analysis of New Malaria Intrahost Models with a Competitive Exclusion Principle\thanks{Received by the editors
October 21, 2005; accepted for publication (in revised form)
July 11, 2006; published electronically December 11, 2006.
 \URL siap/67-1/64327.html}}

\author{Abderrhaman Iggidr\thanks{INRIA-Lorraine and Laboratoire de Math\'ematiques et Applications de Metz UMR CNRS 7122,  University of Metz, 57045 Metz Cedex 01, France (iggidr@math.univ-metz.fr, sallet@loria.fr).} 
\and Jean-Claude Kamgang\thanks{Department of Mathematics, ENSAI,  University of Ngaound\'er\'e, P.O.\ Box 455, Ngaound\'er\'e, Cameroon 
(kamgang@loria.fr).}
\and Gauthier Sallet\footnotemark[2] 
 \and Jean-Jules Tewa\thanks{Department of Mathematics, University of Yaound\'e I, Yaound\'e, Cameroon (tewajules@yahoo.fr).}}

\begin{document}
\slugger{siap}{2006}{67}{1}{260--278}
\maketitle

\setcounter{page}{260}

\begin{abstract}
In this paper we propose a malaria within-host model with $k$ classes of age for the parasitized red blood cells and $n$ strains for the parasite. We provide a global analysis for this model. A competitive exclusion principle holds. If $\mathcal R_0$, the basic reproduction number, satisfies  $\mathcal R_0 \leq 1$, then  the disease-free equilibrium is globally asymptotically stable. On the contrary if $\mathcal R_0 >1$, then generically there is a unique endemic equilibrium which corresponds to the endemic stabilization of the most virulent parasite strain and to the extinction of all the other parasites strains. We prove that this equilibrium is globally asymptotically stable on the positive orthant if a mild sufficient condition is satisfied.
\end{abstract}



\begin{DOI}
10.1137/050643271
\end{DOI}

\pagestyle{myheadings}
\thispagestyle{plain}
\markboth{A. IGGIDR, J.-C. KAMGANG, G. SALLET, AND J.-J. TEWA}{ANALYSIS OF NEW MALARIA INTRAHOST MODELS}

\section{Introduction}
In this paper we consider intrahost models for malaria.  These models describe the interaction  of a parasite, namely a protozoa {\it Plasmodium falciparum}, with its target cells, the red blood cells (RBC). During the past decade there has been considerable work on the mathematical modeling of {\it Plasmodium falciparum} infection 
\cite{9602375,10640432,LLoyd98,Gravenor02JTB,GravenorKwiat98PNAS,GraKwi95,1362993,HetAnd96,11085243,10373354,McKBos97,9735270,11315171,9836304}. A review has been done by Molineaux and Dietz in  \cite{10697860}.  
\enlargethispage{4pt}

We give a brief review of the features of malaria. Malaria in a human begins with an inoculum of {\it Plasmodium} parasites (sporozoites) from a female { \it Anopheles } mosquito. The sporozoites enter the liver within minutes. After a period of asexual reproduction in the liver the parasites (merozoites) are released in the bloodstream where the asexual erythrocyte cycle begins. The merozoites enter RBC, grow, and reproduce over a period of approximately $48$ hours after which the erythrocyte ruptures releasing 8--32 ``merozoites"  daughter parasites that quickly invade  a fresh erythrocyte to renew the cycle. This blood cycle can be repeated many times, in the course of which some of the merozoites instead develop in the sexual form of the parasites: gametocytes. Gametocytes are benign for the host and are waiting for the mosquitoes.

The first mathematical  model of the erythrocyte cycle was proposed by Anderson, May, and Gupta \cite{AnderMay89Parasitology}. This original model has been extended in different directions 
\cite{9602375,AnderMay89Parasitology,LLoyd98,1362993,HetAnd96,11085243,9836304}.

The original model \cite{AnderMay89Parasitology} is given by the following system:\vspace{-6pt}
\begin{equation}\label{AMGori}
 \left \{ \begin{array}{l}
\dot x= \Lambda-\mu_x x - \beta x\,m, \\[-1pt]
\dot y= \beta x \, m -\mu_y \,y,  \\[-1pt]
\dot m=r\,\mu_y \,y -\mu_m \, m -\beta \,x \,m. \\
\end{array}
\right. 
\end{equation}\vspace*{-12pt}

\newpage
\noindent The state variables  are denoted by $x$, $y$, and $m$.  The variable $x$ denotes the concentration of uninfected RBC, $y$ the concentration of parasitized red blood cells 
(PRBC), and $m$ the concentration of the free merozoites in the blood.

%
We briefly sketch  the  interpretation of the parameters. Parameters $\mu_x$, $\mu_y$, and $\mu_m$ are the death rates of the RBC, PRBC, and free merozoites, respectively. The parameter $\beta$ is the contact rate between RBC and merozoites. Uninfected blood cells are recruited at a constant rate $\Lambda$ from the bone marrow and have a natural life-expectancy of 
$\frac{1}{\mu_x}$~days. Death\vspace*{.5pt} of a PRBC results in the release of an average number of $r$ merozoites. Free merozoites die or successfully invade a RBC.

This system  is isomorphic to numerous systems considered in the mathematical modeling of virus dynamics; see \cite{NowAnd2000,PerelsonMatBio93,Perelson99SIAM} and the references therein. Some authors  ignore the loss term $-\beta \,x \,m$
that should appear in the $m$ equation. Indeed without this loss term,  merozoites can infect RBC without themselves being absorbed, and this allows one merozoite {\em to infect more than one RBC}. 

\looseness=-1The original and the derived malaria models were intended to explain observations, namely parasitaemia, i.e., the concentration $y$ of PRBC and also the decrease of the healthy RBC leading to anaemia. An important characteristic of {\it Plasmodium falciparum}, the most virulent malaria parasite, is sequestration. At the halfway point of parasite development, the infected erythrocyte leaves the circulating peripheral blood and binds to the endothelium in the microvasculature of various organs where the cycle is completed. A measurement of {\it Plasmodium falciparum} parasitaemia taken from a blood smear therefore samples young parasites only. Physician treating malaria use the number of parasites in peripheral blood smears as a measure of infection, and this does not give the total parasite burden of the patient. In some respects this is a weak point of the model (\ref{AMGori}). Moreover antimalarial drugs are known to act preferentially on different stages of parasite development. These facts lead some authors to give a general approach to modeling the age structure of {\it Plasmodium} parasites \cite{Gravenor02JTB,GraKwi95,GravenorKwiat98PNAS,15178766}. Their model is a linear catenary compartmental model. This model is based on a finite number of compartments, each representing a stage of development of the parasite inside the PRBC. The models describe only the dynamics of  the morphological stage evolution of the parasites and make no allowance for the dynamics of the healthy RBC. 

In this paper we propose a model which combines the advantages of the two approaches. We also  consider this model with different strains for the parasites. To encompass the different models of the literature we allow, in this model, to ignore or not the loss term in the $m$ equation. To begin we consider the model with one strain:
\begin{equation}\label{combin}
 \left \{ \begin{array}{l}
\dot x= f(x)-\mu_x x- \beta x\,m, \\
\dot y_1= \beta x \, m -\alpha_1\,y_1, \\
\dot y_2=\gamma_1 \, y_1 -\alpha_2\,y_2, \\
\dots\\
\dot y_k= \gamma_{k-1}\,y_{k-1} -\alpha_k \,y_k,\\
\dot m=r\,\gamma_k \,y_k -\mu_m \, m -u \, \beta \,x \,m.\\
\end{array}
\right. 
\end{equation}
In this system $f(x)-\mu_x \, x$ is the density-dependent growth rate of RBC. The other parameters are positive. In the model of Gravenor et al.\ \cite{LLoyd98}  $\alpha_i=\gamma_i+\mu_i$, and hence $\alpha_i > \gamma_i$. We do not need this requirement, which implies that our model is not necessarily a catenary compartmental model. 
In the literature the parameter $u$ takes the values $u=0$ when the loss of the merozoite when it enters a RBC is ignored or takes $u=1$ when this loss is not ignored. In our analysis $u$ is simply a nonnegative parameter. Except   for these generalizations
this system has already been suggested  by Gravenor and Lloyd \cite{LLoyd98} in their reply to the criticism of Saul \cite{9836304}. We provide a global analysis of this system related to the basic reproduction ratio $\mathcal R_0$ of the  considered model. 

One problem is how to decide upon  the number of parasite compartments in the model. A starting point can be the morphological appearance of the parasite. But if the objective is to reflect the distribution of cycle lengths, the number of compartment can be increased to obtain a gamma distribution. Finally the two approaches can be combined: some compartments are for morphological reasons and others are for behavioral reasons. 
Then this  model can also be interpreted as the application of the method of stages (or the linear chain trick) to the life cycle of PRBC \cite{AnderMay89Parasitology,Jacq96,Lloyd01,11589638,11370974,MacDo78}. In other words a chain of compartments is included to generate a distribution of lags. It is also possible to add a class $y_{k+1}$ in order to allow for the production of gametocytes. Different  numbers of stages, ranging from 5 to 48,   are used in \cite{9778631,Gravenor02JTB,GraKwi95,GravenorKwiat98PNAS}. 
\enlargethispage{4pt}

It is well grounded that a  {\it falciparum }
infection consists of distinct parasite genotypes. The model of Anderson, May, and Gupta has been extended in this direction \cite{1362993,Swin96}. With regard  to such features we propose a model with $k$ stages for the infected RBC, production of gametocytes,  and $n$ genotypes, in the population of parasites. 

One of the important  principles of theoretical ecology is the competitive exclusion principle which states that no two species can indefinitely  occupy the same ecological niche 
\cite{Thie89,ButWalt83,MR2000i:92036,MR83d:92066,1362993,LevPim81,6139816,MaySmith74}.  We provide a global analysis of this model and obtain a generic competitive exclusion result within one host individual. This confirms the simulation results obtained in \cite{1362993}.  We compute the basic reproduction ratio $\mathcal R_0$  of the model. For this model there is always a disease-free equilibrium (DFE). To put it more precisely this equilibrium corresponds to the extinction of all  the parasites, 
including the free parasites and the intraerythrocyte parasites. We prove that if $\mathcal R_0 \leq 1$, then the  DFE is globally asymptotically stable (GAS); in other words the parasites are cleared. If $\mathcal R_0 >1$, then, generically, a unique endemic equilibrium exists corresponding to the extinction of all the strains of parasites but one. We prove that this equilibrium is GAS on the positive orthant under  a mild condition.  For example this condition is automatically satisfied when $u=0$ and 
$f(x)=\Lambda-\mu_x\, x$.  When $u \neq 0$ the criteria, obtained for deciding the winning strain,  differs from other results in the literature. To each $i$-strain  can be associated a basic reproduction number $\mathcal R_0^i$ and a threshold $\mathcal T_0^i$. It turns out, when $u\neq 0$, that this is precisely this  threshold $\mathcal T_0^i$ which distinguishes the fate of the strain and  not $\mathcal R_0^i$ at the difference of
\cite{Thie89,MR2000i:92036}.

\looseness=-1The paper is organized as follows. In section~2 we introduce the model with $k$ stages for the infected RBC and one parasite strain,  with and without gametocyte production. We compute the basic reproduction number and provide a stability analysis. 

In section 3 we consider the model of Anderson, May, and Gupta with $n$  distinct genotypes and production of gametocytes. This model with a constant recruitment function for the erythrocytes, two strains, and one class of age  has been proposed in \cite{1362993}. We have studied this model in \cite{Havre05}. Here using the computation of section 2,   we prove for the general $n$ strain $k$ class of age  model  that if $\mathcal R_0 \leq 1$, then the parasites are cleared and if $\mathcal R_0 >1$, then generically the different  genotypes cannot coexist. Namely a unique equilibrium exists, for which only one genotype is positive, and which is GAS on a dense subset  of  the nonnegative orthant. This result confirms the simulations given in \cite{1362993}.

 Global results of stability for the DFE as well  for the endemic equilibrium for epidemic models are not so common 
\cite{MR87c:92046,MR2002c:92034,cit:JaSiKo,0821.92022,SimJacKoo96,Thie83,MR1993355}. Global stability results for the endemic equilibrium have  often been obtained by using  monotone system techniques  \cite{Hir84,LajYo76}. Usually  the Poincar\'e--Bendixson property of monotone systems in dimension $3$ is used \cite{MR2000e:92045,MR95k:92022,MR2001j:92056,0821.92022,MR97e:34093,MR2002m:92032}. Our results generalize the results of \cite{1035.34045}. 
%
%
\section{\boldmath Stability analysis of  a one strain model with $k$ stages}
We consider a general class of systems. The haemopoiesis is a complex system. In the cited references the recruitment of RBC is given by $\Lambda-\mu_x \, x$. In this  paper we will use a more general function 
$\varphi(x)$. In a more complex system the haemopoiesis could be an input coming from another system:
\begin{equation}\label{AMGK}
 \left \{ \begin{array}{l}
\dot x= f(x)-\mu_x x- \beta x\,m =\varphi(x)-\beta \, x \, m,\\
\dot y_1= \beta x \, m -\alpha_1\,y_1, \\
\dot y_2=\gamma_1 \, y_1 -\alpha_2\,y_2, \\
\dots\\
\dot y_k= \gamma_{k-1}\,y_{k-1} -\alpha_k \,y_k,\\
\dot m=r\,\gamma_k \,y_k -\mu_m \, m -u \, \beta \,x \,m.
\end{array}
\right. 
\end{equation}
We denote by $y$ the column vector $(y_1,\ldots,y_k)^T$. The parameter $u$ is nonnegative.  The reason for this parameter is to encompass some malaria models in which the term $-\beta\,x\,m$  can appear or not. In \cite{9602375} Anderson has considered a system without the $\; -\beta \,x \, m$ in the $\dot m$ equation.  In \cite{NowAnd2000} 
all the basic models of virus dynamics are also without this term. One feature of {\it Plasmodium falciparum}, responsible for the deadly case of malaria, is that more than one parasite can invade RBC. In this case $u$ is the mean number of parasites invading RBC and thus disappearing from the circulating blood.

Some authors \cite{1362993,9735270} have included in the model production of gametocytes. In the course of the production of merozoites from bursting erythrocytes, some invading merozoites develop into the sexual, nonreplicating transmission stages known as gametocytes. The gametocytes are benign and transmissible to mosquitoes. We can also, following these authors, include a production of gametocytes in our model. If we denote by $y_{k+1}$ the ``concentration of gametocytes," the model becomes
\begin{equation}\label{AMGKgam}
 \left \{ \begin{array}{l}
\dot x= f(x)-\mu_x x- \beta x\,m =\varphi(x)-\beta \, x \, m,\\
\dot y_1= \beta x \, m -\alpha_1\,y_1, \\
\dot y_2=\gamma_1 \, y_1 -\alpha_2\,y_2, \\
\dots\\
\dot y_k= \gamma_{k-1}\,y_{k-1} -\alpha_k \,y_k,\\
\dot y_{k+1}= \rho \, \gamma_k \,y_k - \alpha_{k+1} \,y_{k+1},\\
\dot m=r\,\gamma_k \,y_k -\mu_m \, m -u \, \beta \,x \,m.
\end{array}
\right. 
\end{equation}
We start to analyze the system with minimal hypothesis on $f$ but nevertheless plausible from the biological point of view. The function $f$ gives  the production of erythrocytes from the bone marrow.
The function $\varphi (x)= f(x)-\mu_x x$ models the population dynamic of RBC in the absence of parasites.  The RBC have a finite lifetime, and  
then $\mu_x$ represents the average per capita death rate of RBC. The function $f$ models  in some way homeostasis. In this paper we suppose that $f$ depends only on $x$. It could be assumed that the recruitment function depends on $x$ and the total population of erythrocytes $x+\sum_i \,y_i$. In this paper we will  analyze  the simplified  case which is the model  considered in all the referenced literature.
The rationale behind  this simplification  is that in a malaria primo-infection typically $y$ is in the order of $10^{-1}$ to $10^{-4}$ of the concentration of healthy erythrocytes $x$. This can be confirmed from the data of malaria therapy. In the last century neurosyphilitic patients  were given malaria therapy, which was routine care at that time. Some of them were infected with {\it Plasmodium falciparum}. Data were collected at the National Institutes of Health laboratories in Columbia, SC and Milledgeville, GA during the period 1940 to 1963 \cite{10432044}.

We assume that $f$ is a $\mathcal C ^1$. Since homeostasis is maintained we assume that the dynamic without parasites is asymptotically stable. In other words, for the system 
$$ \dot x = f(x)-\mu_x \,  x=\varphi (x)$$
there exists a unique $x^* > 0$ such that 
\begin{equation}\label{hypof}
\varphi(x^*)=0,  \;  \; \mbox{ and} \; \;   \varphi(x) >0 \;  \; \mbox{ for} \; \; 0 \leq x < x^*, \;\; \mbox{and} \;\; \varphi(x)  <0 \; \; \mbox{for} \; x > x^*.
\end{equation}
\subsection{Notation}
We will rewrite systems (\ref{AMGK}) and (\ref{AMGKgam}) in a condensed simpler form. 

Before we introduce some classical notation.  

We identify vectors of $\R^n$ with $n \times 1$ column vectors. $\langle \,\mid \, \rangle $ denotes the euclidean inner product. $ \| z \|_2^2= \langle \, z \mid z \,\rangle $ is the usual euclidean norm.

The family $\{e_1,\ldots, e_n\}$ denotes the canonical basis of the vector space $\R^n$. For example $e_1=(1,0, \ldots, 0)^T$. We denote by  $e_{\omega}$ the last vector of the canonical basis, $e_{\omega}=(0,\ldots,0,1)^T$.

If $z \in \R^n$, we denote by $z_i$ the $i$th component of $z$. Equivalently $z_i = \langle \, z \mid e_i \,\rangle$.


For a matrix $A$ we denote by $A(i,j)$ the entry at the  row $i$, column $j$.
For matrices $A,B$ we write $A \leq B$ if $A(i,j) \leq B(i,j)$ for all $i$ and $j$, $A<B$ if $A \leq B$ and $A \neq B$, and $A \ll B$ if $A(i,j) < B(i,j)$ for all $i$ and $j$.

$A^T$ denotes the transpose of $A$. Then $\langle \, z_1\mid z_2 \,\rangle = z_1^T \, z_2$. The notation $A^{-T}$ will denote the transpose of the inverse of $A$.

For this section we rewrite the systems (\ref{AMGK})  and (\ref{AMGKgam}) under a unique form:  
\begin{equation}\label{AMGK2}
 \left \{ \begin{array}{l}
\dot x= \varphi(x)-\beta \, x \, \langle \, e_{\omega} \mid z \,\rangle, \\
\dot z= \beta \, x \, \langle \, e_{\omega} \mid z \,\rangle   \, e_1 +A_0 \, z- u \, \beta \, x \, \langle \, e_{\omega} \mid z \,\rangle   \, e_{\omega}.
\end{array}
\right. 
\end{equation}
In the case of the system (\ref{AMGK}) we have for $A_0$
\begin{equation}\label{A0}
A_0=\begin{bmatrix}
-\alpha_1&0&0&\cdots&0&0 \\
\gamma_1 & -\alpha_2& 0 & \cdots &0&0 \\
0& \gamma_2 & -\alpha_3 & \cdots &0&0\\
\vdots & \ddots & \ddots & \ddots & \vdots&\vdots \\
0 &\cdots & 0 & \gamma_{k-1} & -\alpha_k &0\\
0&\cdots & 0&0&r \gamma_k & \S -\mu_m
\end{bmatrix}
\end{equation}
and an analogous formula for (\ref{AMGKgam}).

We define the matrix $A(x)=A_0-\beta \, x \, e_{\omega} \, e_{\omega}^T$. This a Metzler stable matrix. (A Metzler matrix is a matrix with nonnegative  off-diagonal entries \cite{0815.15016,MR94c:34067,0458.93001}.)

It is not difficult to check that the nonnegative orthant is positively invariant by (\ref{AMGK2}) and that there exists  a compact absorbing set $K$ for this system. An  absorbing set $D$ is a neighborhood such that a trajectory of the system starting from  any initial condition enters and remains in $D$ for  a sufficiently large time $T$.  
\subsection{Global stability results}\label{GASresult}
We can now give the main result of this section. 
%
%
\begin{theorem}\label{stabAMGK}
We consider the system {\rm (\ref{AMGK2})}  with the hypothesis  {\rm (\ref{hypof})} on $\varphi$ satisfied.
We define the basic reproduction ratio of the system {\rm (\ref{AMGK})} and {\rm (\ref{AMGKgam})}  by
\begin{equation}\label{R0}
\mathcal R_0=\frac {  r \beta  x^*} { \mu_m \, +\, u \, \beta\, x^*}\, \frac{\gamma_1 \cdots \gamma_k}{\alpha_1 \cdots \alpha_k}.
\end{equation}\unskip
\begin{remunerate}
\item 
The system {\rm (\ref{AMGK})} is GAS on $\R^{k+2}_+$  
(respectively, {\rm (\ref{AMGKgam})} on $\R^{k+3}_+$) at the DFE $(x^*,0,\ldots,0)$  if and only if $ \mathcal R_0 \leq 1$.
\item If $\mathcal R_0 > 1$, then the DFE is unstable and there exists a unique endemic equilibrium (EE) in the positive orthant,
$(\bar x, \bar z ) \gg 0$, given by
 \begin{equation}\label{endem2}
\left\{
\begin{array}{l}
\bar x = \d{\frac{\mu_m}{\beta \left [ \d{r \frac{\gamma_1 \cdots \gamma_k}{\alpha_1\cdots \alpha_k}       }-u\right ] } ,}\\[8mm]
\bar z =  \varphi (\bar x) \, \left (-A_0\right )^{-1} (e_1-u \, e_{\omega}).
\end{array}
\right.
\end{equation}

Denoting $\alpha^* = -\max_{x\in [0,x^*]} \left ( \varphi'(x) \, \right )$,
if  
\begin{equation}\label{SCstab}
u \, \beta \, \varphi(\bar x) \leq \alpha^* \, \mu_m,
\end{equation}
then  the EE is GAS on the nonnegative orthant, except for initial conditions on the $x$-axis.
\end{remunerate}
\end{theorem}
%

%
\textit{Proof of Theorem}~\ref{stabAMGK}. 
To begin we will consider the system (\ref{AMGK}) without gametocytes, i.e., the system
(\ref{AMGK2}) with $A_0$ as defined in (\ref{A0}). The stability analysis for (\ref{AMGKgam}) follows easily from the stability analysis of  (\ref{AMGK}).

In a first step we will compute  $\mathcal R_0$. 
We use our   preceding notation and define $A^*=A(x^*)$, i.e., the matrix 
computed at the equilibrium $x^*$ of $\varphi$, which is a stable Metzler matrix. We will use, repeatedly in what follows, the property that if $M$ is a stable Metzler matrix, then $-M^{-1} \geq 0$ \cite{0815.15016}. The expression of $\mathcal R_0$ is obtained easily by using the next generation matrix of the system (\ref{AMGK}) \cite{MR1938888,MR1057044,MR2002k:92027}. We have for the basic reproduction number
$$\mathcal R_0= \beta \, x^* \,  \left  \langle \, -\left (A^*\right )^{-1}\, e_1 \mid e_{\omega}  \,\right\rangle. $$
If we remark that the matrix $A^*$ is the matrix $A_0$ modified by a rank-one matrix, namely $A^*=A_0-u \, \beta \, x^* \, e_{\omega}\,e_{\omega}^T$, we can use the Sherman--Morrison--Woodbury formula
$$-\left (A^*\right )^{-1}=-A_0^{-1} - \frac{u \, \beta \, x^*}{1+u \,\beta \, x^* \, e_{\omega}^T \, \left ( -A_0\right )^{-1}\, e_{\omega}} \,  \left ( -A_0\right )^{-1}   \, e_{\omega} \, e_{\omega}^T\,\left (- A_0\right )^{-1}$$
or equivalently
$$-\left (A^*\right )^{-1}=-A_0^{-1} - \frac{ u\,\beta \, x^*}{\mu_m+\beta \, x^* } \,   \, e_{\omega} \, e_{\omega}^T\,\left (- A_0\right )^{-1}.$$
This shows that 
$-\left (A^*\right )^{-1}$ is obtained from $-A_0^{-1}$ by multiplying the last line of $-A_0^{-1}$ by
$\frac{\mu_m}{\mu_m+u \, \beta \, x^*}$. Then  we get
$$\mathcal R_0= \beta \, x^*\, \frac{\mu_m}{\mu_m+u \, \beta \, x^*} \, \left \langle \, -\left (A_0\right )^{-1}\, e_1 \mid e_{\omega}  \,\right \rangle , $$
and then in computing the last entry of the first column of $A_0$ we obtain (\ref{R0}). 

\newpage
We remark  that  $\mathcal R_0 >1$ is equivalent to the following threshold condition:\vspace{4pt} 
\begin{equation}\label{T0}
\mathcal T_0= \frac{\beta \, x^*} {\mu_m} \left [ \mu_m \, \left \langle \, -\left (A_0\right )^{-1}\, e_1 \mid e_{\omega}  \,\right\rangle -u  \right ] =\beta \, x^* \, \left \langle \, -\left (A_0\right )^{-1}\, (e_1-u \, e_{\omega} )\mid e_{\omega} \,\right \rangle>1.
\end{equation}
We are now ready to analyze the stability of the DFE.

It is well known that if $\mathcal R_0 >1$, then the DFE is unstable \cite{MR1057044}, which implies that the condition $\mathcal R_0 \leq 1$ is necessary for stability. 

To prove the sufficiency, in a second  step, we consider  the following function defined on the nonnegative orthant:
\begin{equation}\label{LiapDFE}
 V_{ DFE} (z)= \beta \, x^* \, \langle \, e_{\omega} \mid (-A_0^{-1}) z \,\rangle .
\end{equation}
Its time derivative along the trajectories of system (\ref{AMGK2}) is
\[\dot V_{ DFE} =  \beta \, x \, \langle \, e_{\omega} \mid z \,\rangle \,\beta \, x^*\, \langle \, e_{\omega} \mid (-A_0)^{-1} \,(e_1-u \, e_{\omega}) \,\rangle -\beta \, x^* \langle \, e_{\omega} \mid z\,\rangle
\]
or equivalently, using the expression of $\mathcal T_0$ given in  (\ref{T0}),
\begin{equation}\label{ineqVdotDFE}
\dot V_{ DFE} = \beta \, \langle \, e_{\omega} \mid z \,\rangle \,  \left ( \mathcal T_0 \, x -x^*\right ).
\end{equation}
Now we  take as a candidate Liapunov function, defined on the nonnegative orthant minus the hyperplane face  $x=0$,
\[V=  (x-x^* \ln x) -  x^*(1-\ln x^*) + V_{ DFE} (z).
\]
This function is positive definite (relatively to the DFE) on 
$\R^{k+2}_{+,x>0}=  \{(x,y,m) \in \R^{k+2}_{+}\;:\; x>0\}$.
Its time derivative is given by
\[      \dot V =  \frac{x-x^*}{x} \, \varphi (x) -(x-x^*)\, \beta \, \langle \, e_{\omega} \mid z\,\rangle + \beta \, \langle \, e_{\omega} \mid z \,\rangle \,  \left ( \mathcal T_0 \, x -x^*\right )    \]
or assuming $\mathcal R_0 \leq 1$
$$ \dot V =  \frac{x-x^*}{x} \, \varphi (x) + \beta \, x \, \langle \, e_{\omega} \mid z\,\rangle \, \left ( \mathcal T_0 -1 \right ) \; \leq \;  0.$$
By assumption (\ref{hypof}) we have $(x-x^*) \varphi(x)  \leq 0$ for all $x\geq 0$. Therefore $\dot V \leq 0$ for all 
$(x,z) \in \R^{k+2}_{+,x>0}$, which proves the stability of the DFE. Its attractivity follows from  
LaSalle's invariance principle \cite{Bhatia,Lasalle61,Lasalle76}, since the largest invariant set contained in 
$\{(x,z) \in \R^{k+2}_{+,x>0}\ : \dot V=0 \}$ is reduced to the DFE.
On the other hand the vector field is strictly entrant on the face $x=0$. Hence the whole orthant $\R^{k+2}_{+}$ belongs to the region of attraction of the DFE.
%
%
%

Now we assume that $\mathcal R_0 >1$. The equilibria 
$(\bar x, \bar z)$ of the system, different from the DFE,  are determined by the relations
$$\bar z = \beta \, \bar x \langle \, \bar z \mid e_{\omega} \,\rangle \, \left (- A_0 \right )^{-1} (e_1-u \, e_{\omega}).$$
Replacing $\bar z $ in $\langle \, \bar z \mid e_{\omega} \,\rangle$ we obtain 
\begin{equation}\label{barz}
\langle \, \bar z \mid e_{\omega} \,\rangle = \beta \, \bar x \,  \langle \, \bar z \mid e_{\omega} \,\rangle \,
 \langle \,  \left (- A_0 \right )^{-1} (e_1-u \, e_{\omega}) \mid e_{\omega} \,\rangle.
\end{equation}
If $\langle \, \bar z \mid e_{\omega} \,\rangle=0$, then 
$\varphi (\bar x)=0$, we obtain $\bar x = x^*$, and hence $\bar z = 0$; i.e., the corresponding equilibrium is the DFE.
In the other case, i.e., $\langle \, \bar z \mid e_{\omega} \,\rangle\neq 0$, the relation (\ref{barz}) gives
\begin{equation}\label{barx}
\beta \, \bar x \, \left  \langle \,  \left (- A_0 \right )^{-1} (e_1-u \, e_{\omega}) \mid e_{\omega}  \,\right \rangle=1.
\end{equation}
Using $   \langle \,  \left (- A_0 \right )^{-1} \, e_{\omega} \mid e_{\omega}   \,\rangle = \frac{1}{\mu_m}$ we finally have
$$\bar x = \frac{\mu_m}{\beta \,  \left [  \mu_m \, \langle \, \left (- A_0 \right )^{-1} e_1 \mid e_{\omega}  \,\rangle -u \right ]}= \frac{x^*}{\mathcal T_0}.$$
We deduce that if $\mathcal R_0>1$, then $0<\bar x <x^*$, and hence $\varphi (\bar x)>0$. Therefore 
$$\bar z = \varphi (\bar x) \, \left (-A_0\right )^{-1} (e_1-u \, e_{\omega}).$$
The last component of $\bar z$, 
$\langle \,  \bar z \mid e_{\omega}  \,\rangle= \bar m$, is given by 
$$\bar m=  \frac{\varphi(\bar x)}{\beta \, \bar x} > 0.$$
The $k$ first components of $\bar z$ are given by the $k$ first components of
$\varphi (\bar x) \, \left (-A_0\right )^{-1} e_1$. 
It is straightforward to check that the first column of $(-A_0)^{-1}$ namely 
$(-A_0)^{-1}\,e_1 \gg 0$, which proves that $\bar z \gg 0$.
We have then proved that  there is a unique EE in the positive orthant if and only if $\mathcal R_0 >1$.

Finally we will prove a sufficient condition for the 
global asymptotic stability of the EE. To this end we define the following candidate Liapunov function on the positive orthant minus the face corresponding to $x=0$: 
\begin{equation}\label{LiapEE}
V_{ EE}(x,y,m) =   a(x- \bar x \ln x)+ \sum_{i=1}^k \,  b_i \, (y_{i} -\bar y_{i} \ln y_{i} ) \,+b_{k+1} \, (m-\bar m \, \ln m).
\end{equation}
This function has a unique global minimum in $(\bar x , \bar y , \bar m)$.
We will choose the coefficients $a, b_i,b_{k+1}$ such that in the computation of $\dot V$,   the linear terms in     $y_i$ and $m$ and the bilinear terms in 
$x\,m$  cancel. Let us show that it is possible with positive coefficients.
To this end we rewrite the function $V_{ EE}$ using the notation $z=(y,m)^T$,
$\ln z=(\ln z_1,\ln z_2, \ldots, \ln z_{k+1})^T$, and 
$b=(b_1,\ldots,b_k,b_{k+1})^T$:
\[
V_{ EE}(x,z) =   a(x- \bar x \ln x)+ 
\langle \, b \mid z-diag(\bar z) \ln z\,\rangle. 
\]
Consider the block matrix
$$M=\begin{bmatrix}
-1 & (e_1-u \,e_{\omega})^T \\
\beta \, \bar x \, e_{\omega} & A_0^T
\end{bmatrix}.
$$ 
Using classical Schur complement techniques and the relation (\ref{barx}) on $\bar x$, we have 
$$
\begin{array}{r@{\;}l}
   \text{\rm det}(M)&=\text{\rm det}(A_0)[-1+\beta \, \bar x (e_1-u \, e_{\omega})^T \, (-A_0^{-T})\, e_{\omega}]    \\
  & =\text{\rm det}(A_0)[-1+\beta \, \bar x \,\langle \, -A_0^{-1} \, (e_1-u\,e_{\omega} )\mid e_{\omega} \,\rangle]=0. \\
\end{array}
$$
Since the matrix $M$ is obviously of codimension $1$ ($A_0$ is nonsingular) the kernel of $M$ is of dimension $1$. Then there exists $a \in \R$ and $b \in \R^{k+1}$ such that
\begin{subequations}
\begin{equation} \label{kern1}
a=(e_1-u\,e_{\omega})^T \, b =\langle \, b \mid e_1-u\,e_{\omega}\,\rangle 
\end{equation}
and
\begin{equation}\label{kern2}
b= a \, \beta \, \bar x \,(-A_0^{-T})\,e_{\omega}.
\end{equation}
\end{subequations}
Since the kernel is one dimensional, $a$ can be chosen arbitrarily. Thanks to the structure of $A_0$, if $a >0$, then $b \gg 0$.

The derivative of $V$ along the trajectories of 
 (\ref{AMGK2}) is given by
 \[
\begin{array}{r@{\;}l}
\dot V_{EE} =& \d{
 a \, \frac{x- \bar x}{x}\, \varphi (x) -a \, \beta x \, \langle \, e_{\omega} \mid z \,\rangle + a \, \beta \, \bar x  \, \langle \, e_{\omega} \mid z \,\rangle + \beta \, x \,  \langle \, e_{\omega} \mid z \,\rangle \, \langle \, b \mid e_1-u\,e_{\omega}\,\rangle 
                 }\\[2mm]
& \d{
+\, \langle \, b \mid A_0z \,\rangle + \langle \, b \mid \text{\rm diag} (\bar z) \, \text{\rm diag} (z)^{-1} \, \dot z \,\rangle 
     }\\[2mm]
=&  \d{
 a \, \frac{x- \bar x}{x}\, \varphi (x)
+ \langle \, b \mid \text{\rm diag} (\bar z) \, \text{\rm diag} (z)^{-1} \, \dot z \,\rangle 
     }\\[2mm]
& \d{
+\, a \, \beta \, \bar x  \, \langle \, e_{\omega} \mid z \,\rangle
+ \langle \, b \mid A_0z \,\rangle
+\beta x \, \langle \, e_{\omega} \mid z \,\rangle\,
\bigl(\langle \, b \mid e_1-u\,e_{\omega}\,\rangle  - a \bigr).
       }
\end{array}
\]
Using the relation (\ref{kern2}) we see that
$$\langle \, b \mid A_0\, z \,\rangle=- a \,  \beta \,  \bar x  \, \langle \, (A_0^{-T})\, e_{\omega} \mid A_0 z \,\rangle= - a \, \beta \, \bar x \,  \langle \,  e_{\omega} \mid z \,\rangle.$$
Therefore the linear terms in $z$ cancel. The same is true for the bilinear terms thanks to the relation (\ref{kern1}).
Finally we get
$$ \dot V_{EE} =  a \, \frac{x- \bar x}{x}\, \varphi (x)  + \langle \, b \mid \text{\rm diag} (\bar z) \, \text{\rm diag} (z)^{-1} \, \dot z \,\rangle.$$
We choose $b_{k+1}=1=\langle \, b \mid e_{\omega}\,\rangle=a\, \beta \, \bar x \, \langle \, -A^{-T} \, e_{\omega} \mid e_{\omega}\,\rangle=a \, \beta \, \bar x \frac{1}{\mu_m}$. In other words $a=\frac{\mu_m}{\beta \, \bar x}$. With the hypothesis $\mathcal R_0 >1$ we have $a>0$,  and hence $b \gg 0$ as wanted. With this choice developing $\dot V$ gives
$$
\begin{array}{r@{\;}l}
	 \d \dot V_{ EE} =& \d  a\, f(x) -a \mu_x \,x -a f(x) \frac{\bar x}{x}+ a \, \mu_x \, \bar x-b_1 \, \beta \bar y_1 \frac{xm}{y_1}-\sum_{i=2}^k b_i \gamma_{i-1} y_{i-1} \frac{\bar y_i}{y_i}  \\
&\quad  \d +\, \sum_{i=1}^k b_i \, \alpha_i \, \bar y_i - r \, \gamma_k \,  y_k \ \frac{\bar m}{m}+ u\,  \beta \bar m x + \mu_m \bar m.
\end{array}
$$
We collect some useful relations between our coefficients at the EE. We have from the definitions of $a$ and $b$, since $b_{k+1}=1$, 
\begin{equation} \label{coef1}
\left \{  \begin{array}{l}
a+u=b_1, \\
b_1 \,  \alpha_1= \gamma_1\,  b_2, \\
b_2 \,  \alpha_2= \gamma_2 \, b_3, \\
\cdots  \\
b_{k-1} \, \alpha_{k-1}= \gamma_{k-1}\, b_k, \\
b_k \, \alpha_k =  r\,  \gamma_k.
\end{array} \right.
\end{equation}
From these relations and the properties of the EE $\bar z$  we have
\begin{equation}\label{util1}
b_1 \, \beta \, \bar x \, \bar m=b_i \, \alpha_i \, \bar y_i=b_i \, \gamma_{i-1}  \, \bar y_{i-1}=r \, \gamma_k \, \bar y_k
\end{equation}
and
\begin{equation}\label{util2}
a \alpha_1 \bar y_1=\mu_m\,\bar m.
\end{equation}
Replacing, in the expression of $\dot V$, $a \mu_x \bar x$  by $a \, f(\bar x)-a \beta \bar x \bar m = a f(\bar x) -a \alpha_1 \, \bar y_1$ we obtain 
$$
\begin{array}{r@{\;}l}
\d \dot V_{ EE} =&\d  k r \gamma_k \bar y_k +a f( \bar x)+ a f(x)+ (u\, \beta \bar x \bar m -a \mu_x \bar x )\frac{x}{\bar x}-a f(x) \frac{\bar x}{x} \\
&\quad  \d -\, b_1\, \beta \, \bar x \, \bar m \, \frac{x}{\bar x}\,\frac{m}{\bar m}\,\frac{\bar y_1}{y_1}-\sum_{i=2}^k b_i \gamma_{i-1} \bar y_{i-1} \frac{y_{i-1}}{\bar y_{i-1}} \frac{\bar y_i}{y_i} 
-r \gamma_k \bar y_k \frac{y_k}{\bar y_k} \frac{\bar m}{m}.
\end{array}
$$
Using again the relations between the coefficients we get
$$
\begin{array}{r@{\;}l}
\d \dot V_{ EE} =& \d k r \gamma_k \bar y_k +a f( \bar x)+ a f(x)+ (r \gamma_k \bar y_k -a f( \bar x) )\frac{x}{\bar x}-a f(x) \frac{\bar x}{x} \\[12pt]
&\quad \d -\, r \, \gamma_k\, y_k \, \, \frac{x}{\bar x}\,\frac{m}{\bar m}\,\frac{\bar y_1}{y_1} -\sum_{i=2}^k r \gamma_k \bar y_k \frac{y_{i-1}}{\bar y_{i-1}} \frac{\bar y_i}{y_i} 
-r \gamma_k \bar y_k \frac{y_k}{\bar y_k} \frac{\bar m}{m}
\end{array}
$$
and finally
$$
\begin{array}{r@{\;}l}
\dot V_{ EE} =&\d a \left [ f(x)+f(\bar x )-f(\bar x) \frac{x}{\bar x}-f(x) \frac{\bar x}{x}\right ] 
 \\[12pt]
&\quad  \d +\, r\,\gamma_k\, \bar y_k \left [ k +\frac{x}{\bar x} -\frac{x}{\bar x}\,\frac{m}{\bar m}\,\frac{\bar y_1}{y_1}
-\sum_{i=2}^k \, \frac{y_{i-1}}{\bar y_{i-1}} \frac{\bar y_i}{y_i} -\frac{y_k}{\bar y_k} \frac{\bar m}{m}\right ].
\end{array}
$$
Now we will use the fact that there exists $\xi$ in the open interval  $\xi \in\; ]x,\bar x[$ such that $f(x)=f(\bar x)+(x-\bar x) \, f'(\xi)$. Replacing in the preceding expression gives
$$
\begin{array}{r@{\;}l}
\dot V_{ EE} =&\d a f(\bar x) \left [ 2-\frac{x}{\bar x}- \frac{\bar x}{x}\right ] +a \, f'(\xi)\frac{(x-\bar x)^2}{x} \\[12pt]
&\quad  \d +\, r\,\gamma_k\, \bar y_k \left [ k +\frac{x}{\bar x} -\frac{x}{\bar x}\,\frac{m}{\bar m}\,\frac{\bar y_1}{y_1}
-\sum_{i=2}^k \, \frac{y_{i-1}}{\bar y_{i-1}} \frac{\bar y_i}{y_i} -\frac{y_k}{\bar y_k} \frac{\bar m}{m}\right ].
\end{array}
$$
Using the relations (\ref{coef1})--(\ref{util1}) we have
$$a f(\bar x)= (b_1-u)f(\bar x)=b_1 (\mu_x \bar x + \beta \bar x \bar m) -u\, f(\bar x)=b_1 \, \mu_x \, \bar x+r\gamma_k \bar y_k -u \, f(\bar x).
$$
Replacing in the preceding expression of $\dot V$ gives
$$
\begin{array}{r@{\;}l}
\dot V_{ EE} =&\d \left (b_1\, \mu_x \, \bar x -u \, f(\bar x) \right ) \,  \left [ 2-\frac{x}{\bar x}- \frac{\bar x}{x}\right ] +a \, f'(\xi)\frac{(x-\bar x)^2}{x}  \\[12pt]
&\quad  \d +\, r\,\gamma_k\, \bar y_k \left [ k+2 -\frac{\bar x}{ x} -\frac{x}{\bar x}\,\frac{m}{\bar m}\,\frac{\bar y_1}{y_1}
-\sum_{i=2}^k \, \frac{y_{i-1}}{\bar y_{i-1}} \frac{\bar y_i}{y_i} -\frac{y_k}{\bar y_k} \frac{\bar m}{m}\right ].
\end{array}
$$
This can also be written
\begin{equation}\label{dotVEE}
\begin{array}{r@{\;}l}
\dot V_{ EE}= &\d \Phi(x,y,m)=- \left [ b_1\, \mu_x \, \bar x -u \, f(\bar x) \, - a \, \bar x \, f'(\xi) \right ]\frac{(x-\bar x)^2}{x \bar x} \\[12pt]
& \quad \d +\, r\,\gamma_k\, \bar y_k \left [ k+2 -\frac{\bar x}{ x} -\frac{x}{\bar x}\,\frac{m}{\bar m}\,\frac{\bar y_1}{y_1}
-\sum_{i=2}^k \, \frac{y_{i-1}}{\bar y_{i-1}} \frac{\bar y_i}{y_i} -\frac{y_k}{\bar y_k} \frac{\bar m}{m}\right ].
\end{array}
\end{equation}
The term between brackets in the last expression of $\dot V$ is  
nonpositive by the inequality between the arithmetical mean and the geometrical mean. Then a sufficient condition for $\dot V \leq 0$ is
$$b_1\, \mu_x \, \bar x -u \, f(\bar x) \, -a \, \bar x \, f'(\xi) \geq 0.
$$ 
Moreover with this condition $\dot V$ is negative, except at the EE for the system (\ref{AMGK}). This proves the 
global asymptotic stability of the EE on the positive orthant for the system (\ref{AMGK}). 

The vector field associated with the system is strictly entrant on the faces of the orthant, except the $x$-axis, where it is tangent. The basin of attraction of the EE is then the orthant, except the $x$-axis, which is the stable manifold of the DFE. 

Using the function $\varphi (x)=f(x)-\mu_x\, x$ the preceding condition is equivalent to
$$u \,  \varphi (\bar x) \leq -a \, \bar x \, \varphi' (\xi),$$
or equivalently, replacing $a$ by its value $a=\frac{\mu_m}{\beta \, \bar x}$, the condition becomes
$$u \,  \beta \, \varphi( \bar x) \leq -\mu_m \, \varphi'(\xi).$$
Setting $\alpha^*= - \max_{x \in [0,x^*] } \varphi'(x)$ a sufficient condition for global asymptotic stability of the EE is 
$$\mathcal R_0 > 1 \;\;  \text{ and} \;  \; u \, \beta \, \varphi ( \bar x) \leq  \mu_m \alpha^*.$$
We have proved the theorem for the system without gametocytes. We have seen that $\mathcal R_0$ does not depend on the production of gametocytes. If $\mathcal R_0 \leq 1$,   it is easy,   integrating the linear stable $y_{k+1}$ equations of (\ref{AMGKgam})  from the solutions of (\ref{AMGK}),  to see that the DFE is asymptotically stable and that all the trajectories converge to the equilibrium. The same argument is used when $\mathcal R_0 >1$.
This ends the proof of Theorem~\ref{stabAMGK}.\qquad\endproof 

\begin{remark}
{\rm If this model is a model for  a within-host model of malaria, each  coefficient $\alpha_i$ is made of the mortality of the $i$-class and the rate of transmission in the $i+1$-class: $\alpha_i=\mu_i+\gamma_i$. This implies  that $\gamma_i \leq \alpha_i$. We do not need this assumption, and our conclusions are valid for our more general  model. The only hypothesis is that the parameters of the system are positive.}
\end{remark}

\begin{remark}
{\rm In the proof of Theorem~\ref{stabAMGK} the quantity 
$$\beta \, x^* \, \left \langle \, -\left (A_0\right )^{-1}\, (e_1-u \, e_{\omega} )\mid e_{\omega} \, \right\rangle ,$$ 
which we have called $\mathcal T_0$ when $\mathcal R_0>1$, plays a prominent role. When $\mathcal R_0 \leq 1$ and $u\neq 0$ three cases occur: $0<\mathcal T_0 \leq 1$ or 
 $\mathcal T_0 <0$ or  $\mathcal T_0 = 0$.
 
 In the two first cases we can define $\bar x = \frac{x^*}{\mathcal T_0}$, and we obtain an equilibrium $(\bar x , \bar z)$ of the system which is not in the  nonnegative orthant (either $\bar x <0$ or $\bar z <0$).
 
 In the third case, the computations, done in the proof of Theorem~\ref{stabAMGK}, for the research of an equilibrium show that $\langle \, z \mid e_{\omega} \,\rangle =0$, and hence $z=0$, and finally the equilibrium is the DFE $(x^*, 0)$. 

We introduce a definition of $\mathcal T_0$ that will simplify future computations. The case $\mathcal T_0=0$ is special,  since  $\mathcal T_0=\frac{x^*}{\bar x}$ is no longer true. However this case can be thought, by convention and misuse of language,  as $\bar x = + \infty $.}
\end{remark}

\begin{definition}\label{threshold}
We define for the system {\rm (\ref{AMGK})}  the threshold 
\begin{equation}\label{T_0}
\mathcal T_0= 
\frac{x^{*}}
{\d{
      {\frac{\mu_m}{ \beta \left [r \d{
\frac{\gamma_1 \cdots \gamma_k}
{\alpha_1 \cdots \alpha_k}
								        }
		-u \right ]}
}
 }
 	}=\beta \, x^* \, \left \langle \, -\left (A_0\right )^{-1}\, (e_1-u \, e_{\omega} )\mid e_{\omega}  \,\right\rangle.
\end{equation}
When $\mathcal T_0 \neq 0$ we have also $\mathcal T_0 = \frac{x^{*}} {\bar x}$.
\end{definition}

\begin{remark}
{\rm It should be pointed out that the kind of  Liapunov function defined by (\ref{LiapEE}) has a long history of application to 
Lotka--Volterra models  \cite{Goh76,goh77} and was originally discovered by Volterra himself, although he did not use the vocabulary and the theory of Liapunov functions. Since epidemic models are 
``Lotka--Volterra" like models, the pertinence of this function is not surprising. Similar Liapunov functions have been used in epidemiology \cite{BerCap86,koro,koro04,0778.92020, Webb06},  although with different parameters. We have already used this kind of function in a simplified version of this paper in \cite{Havre05}.}
\end{remark}
\subsection{Comparison with known results}
Our stability result improves the one of De Leenheer and Smith \cite{1035.34045} in two directions:
\begin{remunerate}
  \item We introduce $n$ stages for latent classes.   
\item 
Our sufficient condition for the global asymptotic stability of the endemic equilibrium is weaker than the one provided in \cite{1035.34045}; for instance the sufficient condition given in Theorem~\ref{stabAMGK} is satisfied for malaria parameters given in \cite{AnderMay89Parasitology}, while the condition of \cite{1035.34045} is not satisfied.
\end{remunerate}
\subsection{Application to the original AMG model \cite{AnderMay89Parasitology}}
The original Anderson--May--Guptka model is a three dimensional system 
(\ref{AMGori}) which has the same form as system (\ref{AMGK}) with $f(x)=\Lambda$.
The sufficient condition (\ref{SCstab}) applied to the AMG model (\ref{AMGori}) can be written
\begin{equation}\label{SCstabAMG}
\beta \Lambda \leq \frac{r}{r-1}\, \mu_x\, \mu_m.
\end{equation}
For the system (\ref{AMGori}), it is possible to give a weaker sufficient stability condition.
\begin{proposition}\label{stabAMGori}
If $\mathcal R_0 > 1$ and  $\beta \Lambda \leq (\sqrt{r}+\sqrt{r-1})^2\, \mu_x\, \mu_m$, then the EE is a GAS steady state for system {\rm (\ref{AMGori})} with respect to initial states not on the $x$-axis.
\end{proposition}

Since in general the parameter $r$ is larger than 2 (see, for instance,  \cite{HetAnd96}), we have $(\sqrt{r}+\sqrt{r-1})^2 > \frac{r}{r-1}$.
\begin{proof}
Thanks to the computations done before, we have for system (\ref{AMGori})
$$
\dot V_{ EE}  =
(r-1) \Lambda \left [ 2-\frac{x}{\bar x}- \frac{\bar x}{x}\right ]   
+ r\,\mu_y\, \bar y \left [ 1 +\frac{x}{\bar x}
-\frac{y}{\bar y} \frac{\bar m}{m}
 -\frac{x}{\bar x}\,\frac{m}{\bar m}\,\frac{\bar y}{y}
\right].
$$
Define $X=\frac{x}{\bar x}$ and $S=\frac{y}{\bar y} \frac{\bar m}{m}$. Then one can write
\[\begin{array}{r@{\;}l}
    \dot V_{ EE}  &  =
-(r-1) \Lambda \frac{\left( X-1\right)^2}{X}  
+ r\,\mu_y\, \bar y \left( 1 +X - S -\frac{X}{S}\right)  \\[3mm]
      &   = 
-(r-1) \Lambda \frac{\left( X-1\right)^2}{X}  
+ r\,\mu_y\, \bar y\; \Psi(X,S).
\end{array}  
\]
We have $\Psi(X,S) \geq 0 \Leftrightarrow X\leq S \leq 1 \mbox{ or }
X\geq S \geq 1$. On the other hand 
$\Psi(X,S) \leq \Psi(X,\sqrt{X})=(\sqrt{X}-1)^2$.
Therefore 
\begin{equation}
\label{dotVAMG}
\begin{array}{l}
      \dot V_{EE} \leq (r-1) \Lambda (\sqrt{X}-1)^2
\left(
	\frac{r\,\mu_y\, \bar y}{(r-1) \Lambda} - 
	\left(1+\frac{1}{\sqrt{X}}\right)^2
\right),   \\[3mm]
      \dot V_{ EE} \leq (r-1) \Lambda (\sqrt{X}-1)^2
\left(\sqrt{\frac{r\,\mu_y\, \bar y}{(r-1) \Lambda}} + 1 +\frac{1}{\sqrt{X}}\right)
\left(\sqrt{\frac{r\,\mu_y\, \bar y}{(r-1) \Lambda}} - 1 -\frac{1}{\sqrt{X}}\right).
\end{array}
\end{equation}
We have 
$\frac{\mu_y\, \bar y}{\Lambda}=\frac{\Lambda - \mu_x \bar x}{\Lambda} < 1$.
Hence for 
$X\leq X^*=\frac{x^*}{\bar x}=\frac{(r-1)\beta}{\mu_m}x^*$ we have the following:
$\sqrt{\frac{r\,\mu_y\, \bar y}{(r-1) \Lambda}} - 1 -\frac{1}{\sqrt{X}}
<\sqrt{\frac{r}{(r-1) }} - \frac{\sqrt{\mu_m}}{\sqrt{(r-1)\beta x^*}} -1
\leq 0$, since by assumption $\frac{\beta x^*}{\mu_m}=\frac{\beta\Lambda}{\mu_x\mu_m} \leq 
(\sqrt{r}+\sqrt{r-1})^2$. Therefore, the derivative of $V_{EE}$ along the trajectories of system (\ref{AMGori}) is negative definite on the set 
${\cal D}_0=\{(x,y,m) \in \R^{3}_{+}\;:\; 0<x\leq x^*,\; y>0,\; m>0 \}$. By continuity,
there exists $\epsilon >0$ such that $\dot V_{EE}$ is negative definite on the set ${\cal D}_\epsilon=\{(x,y,m) \in \R^{3}_{+}\;:\; 0<x< x^*+\epsilon,\; y>0,\; m>0 \}$.
The global asymptotic stability of the  EE follows from the fact that ${\cal D}_\epsilon$ is an absorbing set for\break system~(\ref{AMGori}).\qquad\end{proof}

\section{\boldmath The general case: $n$ strains with $k$ classes of parasitized erythrocytes}
We define the following system with $k$ classes and $n$ parasite strains:
\begin{equation}\label{generalKN}
 \left \{ \begin{array}{l}
\displaystyle \dot x= f(x)-\mu_x x- \, x \,  \sum_{i=1}^n  \, \beta_i \,m_i =\varphi(x)- \, x \, \sum_{i=1}^n \, \beta_i \,  \, m_i\\
\text{and  for } i=1, \ldots, n, \\
\dot y_{1,i}= \beta_i x \, m_i -\alpha_{1i}\,y_{1,i}, \\
\dot y_{2,i}=\gamma_{1,i} \, y_{1,i} -\alpha_{2,i}\,y_{2,i}, \\
\dots\\
\dot y_{k,i}= \gamma_{k-1,i}\,y_{k-1,i} -\alpha_{k,i} \,y_{k,i},\\
\dot g_i= \delta_i \, y_{k,i} - \mu_{g_i} \, g_i ,\\
\dot m_i=r_i\,\gamma_{k,i} \,y_{k,i} -\mu_{m_i} \, m_i -u \, \beta_i \,x \,m_i.
\end{array}
\right. 
\end{equation}
%
As in preceding sections we rewrite the system as 
\begin{equation}\label{generalKN2}
 \left \{ \begin{array}{l}
\displaystyle \dot x= \varphi(x)- \, x \, \sum_{i=1}^n \, \beta_i \,  \langle \, z_i \mid  e_{i,\omega}\,\rangle \\
\text{and  for } i=1, \ldots, n, \\
\dot z_{i}= x\,  \beta_i \,  \langle \, z_i \mid  e_{i,\omega}\,\rangle \, e_{i,1} + A_i \, z_i - u \, x\,  \beta_i \,  \langle \, z_i \mid  e_{i,\omega}\,\rangle \, e_{i,\omega},\\
\end{array}
\right. 
\end{equation}
where the matrix $A_{i}$ is the analogous of the matrix $A_0$  defined in section~\ref{GASresult}, but corresponding to  the genotype $i$, and the vectors $e_{i,1}$ and $e_{i, \omega}$ are defined accordingly.  We  drop the index $0$ in $A$ for readability.

\begin{theorem}\label{stabgeneralKN}
We consider the system {\rm (\ref{generalKN})} with the hypotheses  {\rm (\ref{hypof})} satisfied.
We define the basic reproduction ratio $\mathcal R_0$ of the system {\rm (\ref{generalKN})} by
$$\mathcal R_0^i=  \frac {  r_i \beta_i  x^*} { \mu_{m_i} \, + u\, \beta_i\, x^*}\, \frac{\gamma_{1,i} \cdots \gamma_{k,i}}{\alpha_{1,i} \cdots \alpha_{k,i}}$$
and
$$\mathcal R_0 = \max_{i=1,\ldots,n} \, \mathcal R_0^i.$$
\begin{remunerate}
\item 
The system {\rm (\ref{generalKN})} is GAS on $\R^{}_+$  at the DFE $(x^*,0,\ldots,0)$  if and only if $ \mathcal R_0 \leq 1$.
\item If $\mathcal R_0 > 1$, then the DFE is unstable.  If $R_0^i >1$, there exists  an  EE in the nonnegative orthant corresponding to the genotype $i$, the value for the other indexes $j \neq i$ are
$y_j=m_j=0$, and  
 \begin{equation}\label{endemgenKN}
\left\{
\begin{array}{l}
\bar x_i = \d{
\frac{\mu_{m_i}}{\beta_i \left [ r_i 
\d{
\frac{\gamma_{1,i} \cdots \gamma_{k,i}}{\alpha_{1,i}\cdots \alpha_{k,i}}
  }-u\right ] }  
           },\\[26pt]
\bar z_i =  \varphi (\bar x_i) \, \left (-A_{i}\right )^{-1} (e_{i,1}-u \, e_{i,\omega}),\\[6pt]
\bar g_i=\d{
\frac{\delta_i}{\mu_{g_i}}\bar z_{i,k}
           },
\end{array}
\right.
\end{equation}
where we  denote by $\bar z_{i,k}$   the $k${\rm th} component of $\bar z_i$.

\item We assume   $\mathcal R_0 >1$. We define $\mathcal T_0^i $ as in 
	Definition~{\rm \ref{threshold}}.
We assume that the generic conditions  $\mathcal T_0^i \neq \mathcal T_0^j $ are satisfied for $i \neq j$. We suppose that the genotypes have been indexed such that
 $$ \mathcal T_0^1 \,> \mathcal T_0^2 \, \geq   \cdots \geq  \mathcal T_0^n. $$
Then the EE corresponding to $\bar x_1$ is asymptotically stable and the EEs  corresponding to $\bar x_j $  for $j \neq 1$ (for those which  are  in the nonnegative orthant) are unstable.
\item We assume that the preceding hypothesis $\mathcal T_0^1 > \mathcal T_0^j $ is  satisfied with $\mathcal R_0 >1$. We denote it by $\alpha^* = -\max_{x\in [0,x^*]} \left ( \varphi'(x) \, \right )$. Then if
$$ u \, \beta_1 \, \varphi( \bar x_1) \leq \mu_{m_1} \, \alpha^*, $$
the equilibrium $(\bar x_1, \bar y_1,\bar m_1, \bar g_1, 0,\ldots ,0)$ is GAS on the orthant minus the $x$-axis and the faces of the orthant  defined by $y_1=m_1=g_1=0$. In other words  the most virulent strain  is the winner and  the other strains  go extinct.
\end{remunerate}
\end{theorem}
%

\begin{proof}
As in Theorem~\ref{AMGK}  there exists a forward invariant compact absorbing set in the nonnegative orthant for the system (\ref{generalKN}), and hence all the forward  trajectories are bounded. 
The variables $g_i$ do not affect the dynamical evolution of the variables $x,\; y_{i,j},\; m_i$, and so 
we can consider the system without the production of gametocytes. We use the Liapunov function
$$ V_{DFE}(z) = \sum_{i=1}^n \, V_{DFE} (z_i)= \sum_{i=1}^n \,  \beta_i \, x^* \, \langle \, e_{i,\omega} \mid (-A_i^{-1}) z_i \,\rangle. $$ 
Using the system written as (\ref{generalKN2}) and the computation (\ref{ineqVdotDFE}) we easily obtain
$$\dot V_{DFE} = \sum_{i=1}^n \,  \beta_i \, \langle \, e_{i,\omega} \mid z_i \,\rangle \,  \left ( \mathcal T_0^i \, x -x^*\right ).$$
Now we define the Liapunov function on the nonnegative orthant minus the hyperplane face $x=0$
$$V(x,z)=(x-x^*\, \ln x)-x^*(1-\ln x^*)+\sum_{i=1}^n \, V_{DFE} (z_i)$$
which gives
$$
\begin{array}{r@{\;}l}
\dot V  =&\d \frac{x-x^*}{x}\, \varphi(x) +\sum_{i=1}^n \, x^* \beta_i \, \langle \, z_i \mid e_{i, \omega} \,\rangle-
\sum_{i=1}^n \, x \beta_i \, \langle \, z_i \mid e_{i, \omega} \,\rangle \\
& \quad \d +\,\sum_{i=1}^n \,  \beta_i \, \langle \, e_{i,\omega} \mid z_i \,\rangle \,   ( \mathcal T_0^i \, x -x^* ) \\
=& \d  \frac{x-x^*}{x}\, \varphi(x) + \sum_{i=1}^n \,  \beta_i \, \langle \, e_{i,\omega} \mid z_i \,\rangle \,  x \, \left  ( \mathcal T_0^i \, -1 \right ). 
\end{array}
$$
Since $\mathcal R_0^i \leq 1$ for all index $i$, we have $\mathcal T_0^i \leq 1$, and hence 
$\dot V \leq 0$. The conclusion follows by Lasalle's invariance principle and consideration of the boundary of the positive orthant.

Now we assume $\mathcal R_0 >1$. The instability of the DFE follows from the properties of $\mathcal R_0$ \cite{MR1057044}.   We assume that the genotypes are indexed such that their corresponding threshold are in decreasing order $\mathcal T_0^1 > \mathcal T_0^2 \geq \cdots \geq  \mathcal T_0^n$. 

We will define  a Liapunov function  on the nonnegative orthant minus the manifold defined by the equations $x=y_1=m_1=0$. For this we need to  recall the definition of the function $V_{EE}(x,y_1,m_1)$ defined in (\ref{LiapEE}): 
$$V_{ EE}(x,y,m) =   a(x- \bar x \ln x)+ \sum_{i=1}^k \,  b_{1,i} \, (y_{1,i} -\bar y_{1,i} \ln y_{1,i} ) \,+b_{1,k+1} \, (m_1-\bar m_1 \, \ln m_1).$$
The coefficients $(a,b_{1,i})$ are positive and defined from $A_1$ as  in the proof of Theorem~\ref{AMGK}  from section~\ref{GASresult}. We also use the function $V_{EE}$ defined in (\ref{LiapDFE}) to consider 
$$
 V(x,z)=\mathcal T_0^1 \,  V_{EE}(x,z_1)  + a \, \sum_{i=2}^n \, V_{ DFE} (z_i)
$$
or equivalently
$$
 V(x,z)=\mathcal T_0^1 \,  V_{EE}(x,z_1)  + a \, \sum_{i=2}^n \, \beta_i \, x^* \, \langle \, e_{i,\omega} \mid (-A_i^{-1}) \, z_i \,\rangle. 
 $$
Using the relation (\ref{dotVEE}) and (\ref{ineqVdotDFE}), we can compute the derivative of $V$ along the trajectories of (\ref{generalKN2}):
$$\begin{array}{r@{\;}l}
    \dot V =&   \d{  \mathcal T_0^1 \Phi(x,z_1) + \, a \,\mathcal T_0^1\,  \sum_{i=2}^n \, \beta_i \, \bar x_1 \, \langle \, e_{i \omega} \mid z_i \,\rangle  \,-\, a \,\mathcal T_0^1\,  \sum_{i=2}^n \, \beta_i \,  x \, \langle \, e_{i \omega} \mid z_i \,\rangle 
                 }\\[5mm]
      & \quad\d{  +\, a \, \sum_{i=2}^n \, \beta_i \, \langle \, z_i \mid e_{i,\omega} \,\rangle \, \left ( \mathcal T_0^i \,x -x^*\right ).  }
\end{array} 
$$
Using $\mathcal T_0^1 \, \bar x_1 =x^*$ from the Definition~\ref{threshold} for the threshold we get
$$\dot V = \mathcal T_0^1 \Phi(x,z_1) + a \, \sum_{i=2}^n \, \beta_i \, \langle \, z_i \mid e_{i,\omega} \,\rangle \,x \,  \left ( \mathcal T_0^i  -\mathcal T_0^1\right ) \leq 0.$$
By Liapunov theorem this  ends the proof for the stability. The 
global asymptotic stability is obtained by a straightforward use of LaSalle's invariance principle, which ends the proof of Theorem~\ref{stabgeneralKN}.\qquad\end{proof}

\begin{remark}
{\rm In the nongeneric case it can be shown, with the help of  the  Liapunov functions used in the theorem,  that there exists a continuum of stable EE. We omit the proof.

In the generic case, the dynamics of the system are completely determined. The nonnegative orthant is stratified in the union of stable manifolds corresponding to the different equilibria. Only the equilibrium corresponding to the  winning strain has a basin of attraction with a nonempty interior.} 
\end{remark}

\begin{remark}
{\rm We have proved  that the most virulent strain, that is, the strain which maximizes its respective threshold $\mathcal T_0^i$, eliminates the other. We obtain the same kind of result as in \cite{Thie89}, where the  authors consider a $SIR$ model with $n$ strains of parasite. They consider that infection by one parasite strain excludes superinfection by other strains (this is also our case) and induces permanent immunity against all strains in  case of recovery. They also guarantee limited population by considering a recruitment depending on the density in a monotone decreasing way. They find that the strain which maximizes  the basic reproduction ratio eliminates the others.  In the case considered by the authors, actually, using our notation, $\mathcal R_0=\frac{x^*}{\bar x}$. In fact in this model 
$\mathcal T_0$ and $\mathcal R_0$ coincide. This is also the case in our model when $u=0$. Hence our result compares with the result of \cite{Thie89}. However in the case $u\neq 0$ this  is $\mathcal T_0^i$, and not $\mathcal R_0^i$, which distinguishes the fate of the strain. Our result is then different  
from~\cite{Thie89}, where this role is devoted to $\mathcal R_0$. The same kind of remarks apply to \cite{Castillo96} and \cite{MR2000i:92036}.}
\end{remark}

\begin{remark}
{\rm In our model the chains are of equal length for each strain. If the chains are of unequal length, the proof is unchanged. We use equal length for notational convenience. A reason to have unequal length could be to model different behavior for two different strains of the parasite.}
\end{remark}
\section{Conclusion}
In this article we have given a parasitic within-host model and  have provided a stability analysis of this model. 

This model incorporates a number $k$ of compartments for the parasitized target cells and considers $n$ strains for the parasite.  The rationale for including multicompartments can be  multiple. One  reason is to take into account biological reasons, e.g., consideration of  morphological or age classes. The second is for behavioral modeling reasons, 
e.g., to model delays described by gamma distribution functions.

This model has been conceived from malaria infection, since it is well grounded that malaria is a multistrain infection. However other parasitic infections can be considered by this model.

We prove that if the basic reproduction number satisfies $\mathcal R_0 \leq 1$, then the DFE is GAS; i.e., the parasite is cleared from the host. Our stability result when $\mathcal R_0 > 1$ can be summarized as a competitive exclusion principle. To each $i$-strain we associate an individual  threshold condition $ \mathcal T_0^i$  as in 
Definition~\ref{threshold}. If $\mathcal R_0 >1$,  if one strain has its individual  threshold strictly larger than the thresholds of the other strains and if a mild sufficient condition is satisfied (for a constant recruitment,
i.e., $f(x)=\Lambda$, this condition is simply 
$u \beta \Lambda \leq \frac{r}{r-u}\,\mu_x\mu_m$),
 then there exists a GAS equilibrium on the positive orthant. This equilibrium corresponds to the extinction of all  strains, except the strain with the largest threshold.  This winning strain  maximizes the threshold and not its individual  basic reproduction number, which is different from  previous analogous results of the literature.

\end{document}